\documentclass[twocolumn]{autart}

\usepackage{cite}
\usepackage{hyperref}
\usepackage{verbatim}
\usepackage{latexsym}
\usepackage{hyperref}
\usepackage{amsmath}
\usepackage{amsfonts}
\usepackage{amssymb}
\usepackage{graphicx}

\newtheorem{theorem}{Theorem}
\newtheorem{proposition}{Proposition}
\newtheorem{lemma}{Lemma}
\newtheorem{corollary}{Corollary}
\newtheorem{assumption}{Assumption}
\newtheorem{remark}{Remark}

\newcommand{\mb}[1]{\mathbb{#1}}
\newcommand{\mc}[1]{\mathcal{#1}}
\newcommand{\mbf}[1]{\mathbf{#1}}
\newcommand{\bb}{\left(\cdot\right)}
\newcommand{\bbb}{\left(\cdot,\cdot\right)}

\newcommand{\bbbbb}{\left(\cdot,\cdot,\cdot,\cdot\right)}
\newcommand{\N}{\mb{N}}
\newcommand{\R}{\mb{R}}

\newcommand{\Rnm}[2]{\R^{{#1}\times{#2}}}
\newcommand{\support}[2]{\operatorname{h}(#1,#2)}
\newcommand{\supportbb}[1]{\operatorname{h}(#1,\cdot)}

\numberwithin{equation}{section}
\date{\today}

\begin{document}
\begin{frontmatter}
\runtitle{Implicit Rigid Tube MPC}
  \title{The Implicit Rigid Tube Model Predictive Control}
  \author[svr]{Sa\v{s}a~V.~Rakovi\'{c}\thanksref{corr}}  
	\address[svr]{Beijing Institute of Technology, Beijing, China.}
  \thanks[corr]{E--mail: sasa.v.rakovic@gmail.com. Tel.: +44 7799775366.}
\begin{abstract}
A computationally efficient reformulation of the rigid tube model predictive control is developed. A unique feature of the derived formulation is the utilization of the implicit set representations. This novel formulation does not require any set algebraic operations to be performed explicitly, and its implementation requires merely the use of the standard optimization solvers.  
\end{abstract}
\end{frontmatter}
\section{Introduction}
\label{sec:01}

Model predictive control (MPC) is a popular control methodology that has influenced both the theoretical control sphere and the applied control domain~\cite{mayne:rawlings:rao:scokaert:2000,rawlings:mayne:2009,rakovic:2012,mayne:2014,rakovic:2015,rakovic:2019,rakovic:levine:2019}. An area in MPC that has been advanced is robust MPC~\cite{scokaert:mayne:1998,chisci:rossiter:zappa:2001,lofberg:2003,langson:chryssochoos:mayne:rakovic:2004,mayne:seron:rakovic:2005,rakovic:2005,calafiore:campi:2006,goulart:kerrigan:maciejowski:2006,rakovic:kouvaritakis:findeisen:cannon:2011,rakovic:kouvaritakis:cannon:2012,rakovic:kouvaritakis:cannon:panos:findeisen:2012,calafiore:fagaino:2013,rakovic:levine:acikmese:2016,zanon:gros:2021}. A widely accepted robust MPC approach is tube MPC.
Research in tube MPC for linear systems has lead to the development of several of its generations~\cite{langson:chryssochoos:mayne:rakovic:2004,mayne:seron:rakovic:2005,rakovic:2005,rakovic:kouvaritakis:findeisen:cannon:2011,rakovic:kouvaritakis:cannon:2012,rakovic:kouvaritakis:cannon:panos:findeisen:2012,rakovic:levine:acikmese:2016}. Tube MPC methods have been made operational within the context of robust output feedback MPC~\cite{mayne:rakovic:findeisen:allgower:2006,mayne:rakovic:findeisen:allgower:2009}, and their extensions to linear parameter varying and nonlinear systems have also been investigated~\cite{falugi:mayne:2011,mayne:kerrigan:van:falugi:2011,kohler:soloperto:muller:allower:2020,yan:duan:2020,rakovic:dai:xia:2023,ping:yao:ding:li:2021,heydari:farrokhi:2021,hanema:lazar:toth:2020,abbas:mannel:nehoffmann:rostalski:2019}. A beneficial utilization of tube MPC is evident in recent algorithms for robust adaptive MPC~\cite{lorenzen:cannon:allgower:2019,kohler:kotting:soloperto:allgower:muller:2020,zhang:shi:2020} and safe learning~\cite{mckinnon:schoellig:2019,wabersich:zeilinger:2021,gros:zanon:2021}. The practical utility of tube MPC is best evidenced by a variety of effective methods~\cite{farrokhsiar:pavlik:najjaran:2013,prado:torres:cheein:2021,nikou:dimarogonas:2019,oshnoei:kheradmandi:muyeen:hatziargyriou:2021,khaitan:lin:dolan:2021} for mission and motion planning for robotic and multi--agent systems as well as autonomous vehicles. 

Unequivocally, tube MPC has become an active and important research area in its own, which still offers several prominent research challenges including, \emph{inter alia}, computational efficiency and design scalability to higher dimensional problems.  Motivated by these important aspects, we derive herein a computationally efficient reformulation of the rigid tube MPC~\cite{mayne:seron:rakovic:2005}. This novel formulation, referred to as \emph{the implicit rigid tube MPC}, does not require any set algebraic operations appearing in the definitions of the rigid tubes and terminal constraint set to be performed explicitly. The proposed formulation can be implemented by a direct use of the standard LINPROG and QUADPROG functions in MATLAB. 

In our proposal, the set algebraic formalism is used explicitly for the necessary technical analysis, while the underlying representations of the considered sets are deployed in an implicit manner for the computational purpose. Such an approach combines, in an effective manner, best of conceptual and  numerical sides of the set algebraic formalism in order to derive a computationally efficient implicit rigid tube MPC synthesis.   In this context, a novel and unique contribution of our proposal is the use of the implicit rigid state tubes and terminal constraint set. The standard forms of the rigid state tubes and terminal constraint set~\cite{mayne:seron:rakovic:2005} are considered and utilized explicitly for the analysis, while their implicit forms~\cite{rakovic:2022,rakovic:zhang:2023:a,rakovic:zhang:2023:b} are used for the efficient numerical implementation.   The developed implicit rigid tube model predictive control retains all desirable topological and system theoretic properties of its antecedent~\cite{mayne:seron:rakovic:2005}, while its offline and online design stages are  reduced to the standard linear and quadratic programming problems. Such an implementation does not require any of the encountered set algebraic operations to be performed explicitly, and it is computationally efficient and applicable in high dimensions. \emph{To the best of our knowledge, the use of the implicit set representations within the context of tube MPC has not been reported, and, thus, it constitutes a novel,  original and valuable contribution to robust MPC.}

\textbf{Paper Outline.} Section~\ref{sec:02} outlines the key aspects of the ordinary rigid tube MPC~\cite{mayne:seron:rakovic:2005}, and details our objectives. Section~\ref{sec:03} delivers a computationally efficient reformulation of the constraints associated with the implicit rigid tubes. Section~\ref{sec:04} focuses on the implicit rigid tube MPC. Section~\ref{sec:05} provides the related design guidelines in a concise step--by--step manner.  Section~\ref{sec:06} reports relevant numerical experience, illustrates the design on a benchmark control system~\cite{leibfritz:2004}, and delivers conclusions. 

\textbf{Basic Nomenclature and Conventions.}
The sets of real numbers and non--negative integers are denoted by $\R$ and $\N$, respectively. Given $a, b\in\N$ such that $a<b$ we denote $\N_{[a:b]} :=\{a, a + 1, \dots, b-1, b\}$, and we write $\N_b$ for $\N_{[0:b]}$. For a matrix $M\in \Rnm{n}{n}$ and a vector $v\in\R^n$, $M^\top$ and $v^\top$ denote their transposes. The spectral radius $\rho(M)$ of a matrix $M\in \Rnm{n}{n}$ is the largest absolute value of its eigenvalues. A matrix $M\in \Rnm{n}{n}$ is strictly stable if and only if $\rho(M)<1$.  The Minkowski sum of nonempty sets $\mc{X}_j,\ j\in\mc{J}$ in $\R^n$ is
\begin{equation*}
\bigoplus_{j\in\mc{J}} \mc{X}_j:= \{\sum_{j\in\mc{J}}x_j\ :\ \forall j\in\mc{J},\ x_j\in\mc{X}_j\}.
\end{equation*}
The image and preimage of a nonempty set $\mc{X}$ under a matrix of compatible dimensions (or a scalar) $M$ are
\begin{equation*}
M\mc{X}:=\{Mx\ :\ x\in \mc{X}\}\text{ and }M^{-1}\mc{X}:=\{y\ :\ My\in \mc{X}\}.
\end{equation*} 
Likewise, for any integer $k\in\N$,
\begin{equation*}
M^k\mc{X}:=\{M^kx\ :\ x\in \mc{X}\}\text{ and }M^{-k}\mc{X}:=\{y\ :\ M^ky\in \mc{X}\}.
\end{equation*} 
A proper $D$--set in $\R^n$ is a closed convex subset of $\R^n$ that contains the origin in its interior. A proper $C$--set in $\R^n$ is a bounded proper $D$--set in $\R^n$. The intersection of finitely many closed half--spaces is a polyhedral set. A polytopic set is a bounded polyhedral set. 
The support function $\supportbb{\mc{S}}$ of a nonempty, closed, convex subset $\mc{X}$ of $\R^n$ is given, for all $y\in\R^n$, by
\begin{equation*}
\support{\mc{X}}{y}:= \sup_x \{y^\top x\ :\ x\in \mc{X}\}.
\end{equation*} 
The support function of a polyhedral set $\mc{P}=\{p\in\R^n\ :\ \forall i\in\mc{I},\ \varphi_i^\top p\le \mu_i\})$, where $(\varphi_i,\mu_i)\in\R^n\times \R$ for each $i\in\mc{I}$ and $\mc{I}$ is a finite index set, can be evaluated, for any given $y\in\R^n$, by solving the following linear programming problem
\begin{align*}
\text{maximize}\ &y^\top p\\
\text{with respect to}\ &p\\
\text{subject to}\ & \forall i\in\mc{I},\ \varphi_i^\top p\le \mu_i.
\end{align*}
 Finally, we do not distinguish between a variable and its vectorized form in algebraic expressions.

\section{Ordinary Rigid Tube MPC}
\label{sec:02}

\subsection{Regular Linear--Polyhedral--Quadratic Setting}
\label{sec:02.01}
We consider linear discrete time systems 
\begin{equation}
\label{eq:02.01}
  x^+=Ax+Bu+w,
\end{equation}
where $x\in \R^n$, $u\in \R^m$ and $w\in \R^n$ are, respectively, the current state, control and disturbance, and $x^+ \in\R^n$ is the successor state. 
\begin{assumption}
\label{ass:02.01}
The matrix pair $(A,B)\in \Rnm{n}{n}\times\Rnm{n}{m}$ is known exactly, and it is strictly stabilizable.
\end{assumption}

We consider general polyhedral stage constraints on the state and control variables $x$ and $u$, as specified by
\begin{equation}
\label{eq:02.02}
(x,u)\in \mc{Y}:=\{(x,u)\ :\ \forall i\in\mc{I}_\mc{Y},\ c_i^\top x+d_i^\top u\le 1\}.
\end{equation}

\begin{assumption}
\label{ass:02.02}
The representation of the polyhedral proper $D$--set $\mc{Y}$ in~\eqref{eq:02.02} is irreducible. The index set $\mc{I}_\mc{Y}:=\{1,2,\ldots,p\}$ is finite. Each pair $(c_i,d_i)\in\R^n\times\R^m,\ i\in\mc{I}_\mc{Y}$ is known exactly.
\end{assumption}

We consider geometric stage bounds on the disturbance variable $w$ induced by a polytopic proper $C$--set $\mc{W}$, i.e.,
\begin{equation}
\label{eq:02.03}
w\in \mc{W}:=\{w\in\R^n\ :\ \forall i\in\mc{I}_\mc{W},\ e_i^\top w\le 1\}.
\end{equation}

\begin{assumption}
\label{ass:02.03}
The representation of the polytopic proper $C$--set $\mc{W}$ in~\eqref{eq:02.03} is irreducible. The index set $\mc{I}_\mc{W}:=\{1,2,\ldots,q\}$ is finite. Each $e_i\in\R^n,\ i\in\mc{I}_\mc{W}$ is known exactly.
\end{assumption}

In terms of the system description, at any time instance $k\in\N$, the system~\eqref{eq:02.01} reads as $x_{k+1}=Ax_k+Bu_k+w_k$. Likewise, in terms of  the available information for control synthesis, at any time instance $k\in\N$, when the decision on the respective current control action $u_k$ is taken, the respective current state $x_k$ is known exactly, while the respective current disturbance $w_k$ and respective future disturbances $w_{k+i},\ i\in\N,\ i\ge 1$ are not known but are guaranteed to obey the set--membership relation $w\in\mc{W}$, i.e., it is guaranteed that, for all $i\in\N$, $w_{k+i}\in\mc{W}$. 

We also consider the quadratic and jointly strictly convex stage cost function $\ell\bbb$ specified, for all $z\in \R^n$ and all $v\in \R^m$, by
\begin{equation}
\label{eq:02.04}
\ell(z,v)=z^\top Qz+v^\top Rv.
\end{equation}

\begin{assumption}
\label{ass:02.04}
The matrices $Q\in \Rnm{n}{n}$ and $R\in \Rnm{m}{m}$ are known exactly, symmetric and positive definite, i.e., $Q=Q^\top \succ 0$ and $R=R^\top \succ 0$.
\end{assumption}

\subsection{Ordinary Rigid Tube Optimal Control}
\label{sec:02.02}
The rigid tube MPC is best seen as a repetitive decision making process, in which the basic decision making process is the rigid tube optimal control. In the rigid tube optimal control, the rigid state tubes and associated affine control policy are employed for the predictions under uncertainty. The rigid state tube $\mbf{X}_N$ is a set sequence $\{\mc{X}_k\}_{k=0}^N$, each term of which is parametrized via a point $z_k\in\R^n$ and an \emph{a priori} specified set $\mc{S}\subset \R^n$ as follows
\begin{equation}
\label{eq:02.05}
\forall k\in\N_N,\ \mc{X}_k:=z_k\oplus \mc{S}.
\end{equation} 
An affine control policy $\Pi_{N-1}\bb$, associated with the rigid state tube $\mbf{X}_N$ and a state $x\in\R^n$, is  a sequence of control laws $\{\pi_k\bbbbb\}_{k=0}^{N-1}$, each of which is parametrized via points $x_k\in\R^n$, $z_k\in\R^n$ and $v_k\in\R^m$ and an \emph{a priori} specified control matrix $K_\mc{S}\in\Rnm{m}{n}$ as follows
\begin{align}
&\forall k\in\N_{N-1},\  \forall x_k\in \mc{X}_k,\nonumber\\
\label{eq:02.06}
&\pi_k(x_k,z_k,v_k,x):=v_k+K_\mc{S}(x_k-z_k).
\end{align} 
The sequences $\{z_k\}_{k=0}^N$ and $\{v_k\}_{k=0}^{N-1}$ are traditionally referred to as the sequences of the rigid state and control tubes centres~\cite{mayne:seron:rakovic:2005}.  The matrix $K_\mc{S}$ and the set $\mc{S}$ are subject to mild conditions detailed in the next subsection. 

For a given $x\in\R^n$, the ordinary rigid tube optimal control takes the form of the following optimization problem
\begin{align}
\text{minimize}\ &\sum_{k=0}^{N-1} (z_k^\top Q z_k+v_k^\top R v_k) + z_N^\top P z_N\nonumber\\
\text{with respect to}\ &\{z_k\}_{k=0}^N\text{ and }\{v_k\}_{k=0}^{N-1}\nonumber\\
\text{subject to}\ &x\in z_0\oplus \mc{S},\nonumber\\
&\forall k\in\N_{N-1},\ z_{k+1}=Az_k+Bv_k,\nonumber\\
&\forall k\in\N_{N-1},\ (z_k,v_k)\in\mc{Y}_{\mc{S}},\text{ and}\nonumber\\
\label{eq:02.07}
&z_N\in\mc{Z}_f.
\end{align}
In this formulation of the rigid tube optimal control, the set $\mc{Y}_\mc{S}$ is a modified stage constraint set, the set $\mc{Z}_f$ is a terminal constraint set and the function $V_f\bb$ specified, for all $z\in\R^n$, by 
\begin{equation}
\label{eq:02.08}
V_f(z)=z^\top P z
\end{equation}
is a terminal cost function. These additional design ingredients are also subject to mild conditions detailed in the next subsection. In the ordinary rigid tube MPC~\cite{mayne:seron:rakovic:2005}, under mild conditions that are outlined in the next subsection, the ordinary rigid tube optimal control problem~\eqref{eq:02.07} reduces to a strictly convex quadratic programming problem and the resulting ordinary rigid tube MPC ensures the highly desirable robust positive invariance and robust exponential stability properties. 

\subsection{Standard Ingredients in Ordinary Rigid Tube MPC}
\label{sec:02.03}

The basic assumption on the matrix $K_\mc{S}$ is the following.
\begin{assumption} 
\label{ass:02.05}
 The matrix $K_\mc{S}\in\Rnm{m}{n}$ is known exactly and such that the matrix $A+BK_\mc{S}$ is strictly stable.
\end{assumption}

The standard choice for the set $\mc{S}$ is a robust positively invariant approximation~\cite{rakovic:kerrigan:kouramas:mayne:2004b} of the minimal robust positively invariant set for the local uncertain $s$--dynamics
\begin{equation}
\label{eq:02.09}
s^+=(A+BK_\mc{S})s+w,
\end{equation} 
which is obtained from the uncertain $s$ system $s^+=As+B\nu+w$ under the local linear feedback control law $\nu=K_\mc{S}s$. Such a choice is captured by the following.
\begin{assumption}
\label{ass:02.06}
 The set $\mc{S}$ is given by
 \begin{equation}
 \label{eq:02.10}
 \mc{S}= (1-\alpha)^{-1}\bigoplus_{j=0}^{N_\mc{S}-1} (A+BK_\mc{S})^j\mc{W},
 \end{equation}
 where a scalar $\alpha \in [0,1)$ and a finite integer $N_\mc{S}\in\N$ satisfy
 \begin{equation}
 \label{eq:02.11}
 (A+BK_\mc{S})^{N_\mc{S}}\mc{W}\subseteq \alpha \mc{W}.
 \end{equation} 
\end{assumption}
The set $\mc{S}$ is also required to be strictly admissible with respect to stage constraints $\mc{Y}$, i.e., that, for all $s\in\mc{S}$, we have $(s,K_\mc{S}s)\in \operatorname{interior}(\mc{Y})$, or, equivalently, that, for all $i\in \mc{I}_\mc{Y}$, for all $s\in \mc{S}$, we have $(c_i^\top +d_i^\top K_\mc{S})s<1$. This can be equivalently stated in a compact form via~\eqref{eq:02.12}.
\begin{assumption}
\label{ass:02.07}
 The set $\mc{S}$ is such that
 \begin{equation}
 \label{eq:02.12}
 \forall i\in \mc{I}_\mc{Y},\  \support{\mc{S}}{c_i+K_\mc{S}^\top d_i}<1.
 \end{equation} 
\end{assumption}

The modified stage constraint set $\mc{Y}_\mc{S}$ is defined by
\begin{equation}
\label{eq:02.13}
\mc{Y}_{\mc{S}}:=\{(z,v)\ :\ \forall s\in\mc{S},\  (z+s,v+K_\mc{S}s)\in\mc{Y}\}.
\end{equation}

The terminal constraint set $\mc{Z}_f$ is the maximal positively invariant set for the virtual terminal $z$--dynamics
\begin{equation}
\label{eq:02.14}
z^+ =(A+BK_\mc{Z})z,
\end{equation}
which is obtained from the disturbance free $z$ system $z^+=Az+Bv$ under the terminal linear feedback control law $v=K_\mc{Z}z$, and whose $z$ is subject to constraints
$(z,K_\mc{Z}z)\in\mc{Y}_\mc{S}$ that can be expressed more compactly  as
\begin{equation}
\label{eq:02.15}
z\in\mc{Z}_\mc{S}=\{z\ :\ (z,K_\mc{Z}z)\in\mc{Y}_\mc{S}\}.
\end{equation}
The basic assumption on the matrix $K_\mc{Z}$ is the following.
\begin{assumption} 
\label{ass:02.08}
 The matrix $K_\mc{Z}\in\Rnm{m}{n}$ is known exactly and such that the matrix $A+BK_\mc{Z}$ is strictly stable.
\end{assumption}

The structural form of the utilized set $\mc{Z}_f$ is well understood~\cite{gilbert:tan:1991,rakovic:zhang:2023:a,rakovic:zhang:2023:b}, and it is captured by the following.
\begin{assumption}
\label{ass:02.09}
 The set $\mc{Z}_f$ is given by
 \begin{equation}
 \label{eq:02.16}
 \mc{Z}_f:=\bigcap_{j=0}^{N_\mc{Z}}(A+BK_\mc{Z})^{-j}\mc{Z}_\mc{S},
 \end{equation}
where a finite integer $N_\mc{Z}\in\N$ satisfies
 \begin{equation}
 \label{eq:02.17}
 \bigcap_{j=0}^{N_\mc{Z}}(A+BK_\mc{Z})^{-j}\mc{Z}_\mc{S}\subseteq (A+BK_\mc{Z})^{-(N_\mc{Z}+1)}\mc{Z}_\mc{S}.
 \end{equation} 
\end{assumption}
Finally, the terminal cost function $V_f\bb$ is required to be a Lyapunov function for the virtual terminal $z$--dynamics $z^+=(A+BK_\mc{Z})z$ over the terminal constraint set $\mc{Z}_f$, which satisfies the usual decrease condition induced by the stage cost function. This is ensured by the following.
\begin{assumption} 
\label{ass:02.10}
 The matrix $P$ is known exactly, symmetric and positive definite (i.e., $P=P^\top \succ 0$), and it is such that
 \begin{equation}
 \label{eq:02.18}
 (A+BK_\mc{Z})^\top P (A+BK_\mc{Z}) - P \preceq -(Q+K_\mc{Z}^\top R K_\mc{Z}).
 \end{equation}
\end{assumption}

\subsection{Article Objectives}
\label{sec:02.04}

The use of the postulated sets $\mc{S}$ and $\mc{Z}_f$ is motivated by a natural desire to reduce the spread of the uncertainty and enlarge the related domain of the attraction. In the current state of the affairs, a common belief has emerged that the explicit representations of the sets $\mc{S}$ and $\mc{Z}_f$ need to be constructed and also utilized for the design of the ordinary rigid tube MPC~\cite{mayne:seron:rakovic:2005}. Inevitably, the actual computation of the sets $\mc{S}$ and $\mc{Z}_f$ and their explicit representations would certainly render the ordinary rigid tube MPC~\cite{mayne:seron:rakovic:2005} computationally impracticable for anything but geometrically rather simple problems. 

Our main conceptual concern is to dispel such an untenable meme impairing the actual utility of the rigid tube MPC. Our main technical objective is to develop a computationally efficient reformulation of the ordinary rigid tube MPC~\cite{mayne:seron:rakovic:2005}, which can be implemented via standard convex optimization solvers without the need to perform explicitly any of the underlying set algebraic operations.

\section{Implicit Rigid Tubes: Constraints Handling}
\label{sec:03}

With the use of the implicit set representations, the constraints appearing in the ordinary rigid tube MPC~\eqref{eq:02.07} can be reformulated as a tractable set of affine equalities and inequalities in a manner that is computationally efficient and plausibly scalable to high dimensions. 

\subsection{Initialization}
\label{sec:03.01}
In light of the utilized form~\eqref{eq:02.10} of the set $\mc{S}$, the rigid state tube initialization constraints, specified via constraints $x\in z_0\oplus \mc{S}$ in~\eqref{eq:02.07},  read as
\begin{equation}
\label{eq:03.01}
x\in z_0\oplus (1-\alpha)^{-1}\bigoplus_{j=0}^{N_\mc{S}-1} (A+BK_\mc{S})^j\mc{W},
\end{equation}
which is a set--membership condition with respect to the Minkowski sum of the finitely many set images. Such a set--membership condition is equivalent to a computationally convenient existence problem, and it does not require the considered set to be constructed explicitly. 
\begin{lemma}
\label{lemma:03.01}
Let $\mc{X}$ and $M$ be a proper $C$--set in $\R^n$ and a matrix in $\Rnm{n}{n}$, respectively. Let also $\mc{J}:=\{0,1,\ldots,j\}$ for a finite  $j\in\N,\ j>0$. For all $y\in\R^n$,
\begin{equation*}
y\in \bigoplus_{j\in\mc{J}} M^j\mc{X}
\end{equation*}
if and only if, for a collection of points $\{x_j\ :\ j\in\mc{J}\}$, 
\begin{equation*}
y=\sum_{j\in\mc{J}} M^jx_j\ \text{with}\  \forall j\in\mc{J},\  x_j\in\mc{X}.
\end{equation*}
\end{lemma}
Within our setting, due to Lemma~\ref{lemma:03.01}, the set--membership condition~\eqref{eq:03.01} holds true if and only if
\begin{align}
&x= z_0+ (1-\alpha)^{-1}\sum_{j=0}^{N_\mc{S}-1} (A+BK_\mc{S})^j\omega_j\  \text{with}\nonumber\\
\label{eq:03.02}
&\forall j\in\N_{N_\mc{S}-1},\  \forall i\in\mc{I}_\mc{W},\  e_i^T\omega_j\le 1.
\end{align}
The direct implementation of the derived equivalent reformulation~\eqref{eq:03.02} requires the consideration of $N_\mc{S}$ decision variables $\omega_j\in\R^n,\ j\in\N_{N_\mc{S}-1}$ and the use of $n$ affine equalities and $N_\mc{S} q$ affine inequalities. 

\subsection{Guaranteed Dynamical Consistency}
\label{sec:03.03}
By~\cite[Theorem~1]{rakovic:kerrigan:kouramas:mayne:2004b}, the utilized set $\mc{S}$ is a polytopic proper $C$--set in $\R^n$ such that $(A+BK_\mc{S})\mc{S}\oplus \mc{W}\subseteq \mc{S}$. Thus, as in the ordinary rigid tube MPC~\cite{mayne:seron:rakovic:2005},  the guaranteed rigid state tube dynamical consistency is simply specified via constraints, for all $k\in\N_{N-1}$, $z_{k+1}=Az_k+Bv_k$ in~\eqref{eq:02.07}, i.e.,
\begin{equation}
\label{eq:03.03}
\forall k\in\N_{N-1},\  z_{k+1}=Az_k+Bv_k, 
\end{equation}
which is $Nn$ affine equalities.
 
\subsection{Stage Constraints}
\label{sec:03.04}

The robust stage constraints satisfaction specified via constraints, for all $k\in\N_{N-1}$, $(z_k,v_k)\in\mc{Y}_\mc{S}$ in~\eqref{eq:02.07} and the definition of the set $\mc{Y}_\mc{S}$ in~\eqref{eq:02.13} require the evaluation of the support function of the Minkowski sum of the finitely many set images. Such an evaluation of the support function can be done efficiently, and it does not require the considered set to be constructed explicitly.
\begin{lemma}
\label{lemma:03.02}
Let $\mc{X}$ and $M$ be a proper $C$--set in $\R^n$ and a matrix in $\Rnm{n}{n}$, respectively. Let also $\mc{J}:=\{0,1,\ldots,j\}$ for a finite  $j\in\N,\ j>0$.  For all $y\in\R^n$,
\begin{equation*}
\support{\bigoplus_{j\in\mc{J}} M^j\mc{X}}{y}=\sum_{j\in\mc{J}}\support{\mc{X}}{(M^j)^\top y}.
\end{equation*}
\end{lemma}
In light of Lemma~\ref{lemma:03.02}, the topological and structural properties of the set $\mc{Y}_\mc{S}$ can be summarized by the following. 
\begin{proposition}
\label{prop:03.01}
Suppose Assumptions~\ref{ass:02.01}--\ref{ass:02.07} hold. The set $\mc{Y}_\mc{S}$ specified by~\eqref{eq:02.13} is a polyhedral proper $D$--set in $\R^{n+m}$ with a, possibly redundant, representation
\begin{equation}
\label{eq:03.04}
\mc{Y}_\mc{S}=\{(z,v)\ :\ \forall i\in\mc{I}_\mc{Y},\  c_i^\top z + d_i^\top v \le 1-f_i\},
\end{equation}
where, for all $i\in\mc{I}_\mc{Y}$, with $\eta_i:=c_i+K_\mc{S}^\top d_i$, the scalars
\begin{align}
f_i:=&\support{\mc{S}}{\eta_i}
\nonumber\\
\label{eq:03.05}
=&(1-\alpha)^{-1}\sum_{j=0}^{N_\mc{S}-1}\support{\mc{W}}{((A+BK_\mc{S})^j)^\top \eta_i}
\end{align}
are such that $f_i\in [0,1)$ (so that $(1-f_i)\in (0,1]$).
\end{proposition}
The representation of the set $\mc{Y}_\mc{S}$ can be constructed efficiently since, in light of Lemma~\ref{lemma:03.02}, for any $\eta\in\R^n$, the evaluation of the value  $\support{\mc{S}}{\eta}$ of the support function $\supportbb{\mc{S}}$  requires merely the evaluation of the values $\support{\mc{W}}{((A+BK_\mc{S})^j)^\top \eta},\ j\in\N_{N_\mc{S}-1}$ of the support function $\supportbb{\mc{W}}$, as utilized in~\eqref{eq:03.05}. 
Due to the preceding construction, the robust stage constraints satisfaction, specified via constraints, for all $k\in\N_{N-1}$, $(z_k,v_k)\in\mc{Y}_\mc{S}$ in~\eqref{eq:02.07}, equivalently reads as
\begin{equation}
\label{eq:03.06}
\forall k\in\N_{N-1},\  \forall i\in\mc{I}_\mc{Y},\  c_i^Tz_k+d_i^Tv_k\le 1-f_i,
\end{equation}
which is $Np$ affine inequalities.

\subsection{Terminal Constraints}
\label{sec:03.05}
The robust terminal constraints satisfaction specified via constraints $z_N\in\mc{Z}_f$ in~\eqref{eq:02.07}, with the form of the set $\mc{Z}_f$ in~\eqref{eq:02.16} in mind, is a set--membership condition with respect to the set intersection of the finitely many set preimages. Such a set--membership condition can be converted into an equivalent, and computationally convenient, existence problem, and it does not require the considered set to be constructed explicitly.
\begin{lemma}
\label{lemma:03.03}
Let $\mc{X}$ and $M$ be a proper $D$--set in $\R^n$ and a matrix in $\Rnm{n}{n}$, respectively. Let also $\mc{J}:=\{0,1,\ldots,j\}$ for a finite $j\in\N,\ j>0$.  For all $y\in\R^n$, 
\begin{equation*}
y\in \bigcap_{j\in\mc{J}} M^{-j}\mc{X}
\end{equation*}
if and only if, for a collection of points $\{y_j\ :\ j\in\mc{J}\}$, 
\begin{equation*}
\forall j\in\mc{J},\  y_j\in\mc{X}\ \text{with}\  y_j=M^jy.
\end{equation*}
\end{lemma}
To make use of Lemma~\ref{lemma:03.03}, we note that, due to form of the set $\mc{Y}_\mc{S}$ obtained in Proposition~\ref{prop:03.01}, the set $\mc{Z}_\mc{S}$ in~\eqref{eq:02.15} is guaranteed to be at least a polyhedral proper $D$--set in  $\R^n$ with a, possibly redundant, representation
\begin{equation}
\label{eq:03.07}
\mc{Z}_\mc{S}=\{z\ :\ \forall i\in\mc{I}_\mc{Y},\  (c_i^\top + d_i^\top K_\mc{Z})z \le 1-f_i\}.
\end{equation}  
Thus, by Lemma~\ref{lemma:03.03}, the set--membership condition
\begin{equation}
\label{eq:03.08}
z_N\in \mc{Z}_f=\bigcap_{j=0}^{N_\mc{Z}}(A+BK_\mc{Z})^{-j}\mc{Z}_\mc{S}
\end{equation} 
holds true if and only if
\begin{align}
&\forall k\in\N_{N_\mc{Z}-1},\  z_{N+k+1}=(A+BK_\mc{Z})z_{N+k}\  \text{with}\nonumber\\
\label{eq:03.09}
&\forall k\in\N_{N_\mc{Z}},\  \forall i\in\mc{I}_\mc{Y},\  (c_i^T+d_i^TK_\mc{Z})z_{N+k}\le 1-f_i.
\end{align}
The direct implementation of the derived equivalent reformulation~\eqref{eq:03.09} requires the consideration of $N_\mc{Z}+1$ decision variables $z_{N+k}\in\R^n,\ k\in\N_{N_\mc{Z}}$ and the use of $N_\mc{Z}n$ affine equalities and $(N_\mc{Z}+1)p$ affine inequalities.  

\section{Implicit Rigid Tube MPC}
\label{sec:04}

We now detail  a computationally efficient implementation of the implicit rigid tube MPC.  \emph{In what follows, for typographical convenience, we deploy the symbol $T$ to denote a triplet $(N,N_\mc{S},N_\mc{Z})$, i.e., $T:=(N,N_\mc{S},N_\mc{Z})$.}

\subsection{Implicit Rigid Tube Optimal Control}
\label{sec:04.01}
The overall decision variable $\mbf{d}_T$ consists of the ``usual" sequences $\mbf{z}_{N-1}:=\{z_k\in\R^n\}_{k=0}^{N-1}$ and $\mbf{v}_{N-1}:=\{v_k\in\R^m\}_{k=0}^{N-1}$, and the ``additional" sequences $\mbf{z}_{[N:N+N_\mc{Z}]}:=\{z_k\in\R^n\}_{k=N}^{N+N_\mc{Z}}$ and $\boldsymbol{\omega}_{N_\mc{S}-1}:=\{\omega_k\in\R^n\}_{k=0}^{N_\mc{S}-1}$. Thus, the decision variable $\mbf{d}_T$ and its dimension $n_{\mbf{d}_T}$ satisfy 
\begin{align}
\mbf{d}_T:&=(\mbf{z}_{N-1},\mbf{v}_{N-1},\mbf{z}_{[N:N+N_\mc{Z}]},\boldsymbol{\omega}_{N_\mc{S}-1})\text{ and}\nonumber\\
\label{eq:04.01}
n_{\mbf{d}_T}&=Nn+Nm+(N_\mc{Z}+1)n+N_\mc{S}n.
\end{align}
The overall cost function $\mbf{d}_T\mapsto V_T(\mbf{d}_T)$ is given by
\begin{align}
V_T(\mbf{d}_T)=&\sum_{k=0}^{N-1} (z_k^\top Q z_k+v_k^\top R v_k)+\sum_{k=N}^{N+N_{\mc{Z}}-1} z_k^\top Q_\mc{Z} z_k\nonumber\\
\label{eq:04.02}
&+ z_{N+N_{\mc{Z}}}^\top P z_{N+N_{\mc{Z}}},
\end{align}
where $Q_\mc{Z}:=Q+K_\mc{Z}^\top R K_\mc{Z}$ is such that $Q_\mc{Z}=Q_\mc{Z}^\top \succ 0$.
\begin{remark}
\label{rem:04.01}
When $N_\mc{Z}=0$, the terminal constraint set $\mc{Z}_f$ is $\mc{Z}_\mc{S}$ and the related terminal constraints~\eqref{eq:03.09} simplify to $z_N\in\mc{Z}_\mc{S}$.  Furthermore, in such a case, $\mbf{z}_{[N:N+N_\mc{Z}]}=\{z_N\}$, and, by convention, $\N_{N_{\mc{Z}}-1}=\emptyset$ and the symbol $\sum_{k=N}^{N+N_{\mc{Z}}-1}$ denotes the sum that has no terms. 
\end{remark}

In our setting, for a given $x\in\R^n$, the implicit rigid tube optimal control problem $\mathfrak{P}_T(x)$ takes the following form
\begin{align}
&\text{minimize}\ V_T(\mbf{d}_T)\nonumber\\
&\text{with respect to}\ \mbf{d}_T=(\mbf{z}_{N-1},\mbf{v}_{N-1},\mbf{z}_{[N:N+N_\mc{Z}]},\boldsymbol{\omega}_{N_\mc{S}-1})\nonumber\\
&\text{subject to}\nonumber\\
&x= z_0+ (1-\alpha)^{-1}\sum_{j=0}^{N_\mc{S}-1} (A+BK_\mc{S})^j\omega_j,\nonumber\\
&\forall j\in\N_{N_\mc{S}-1},\quad  \forall i\in\mc{I}_\mc{W},\quad e_i^T\omega_j\le 1,\nonumber\\
&\forall k\in\N_{N-1},\quad z_{k+1}=Az_k+Bv_k,\nonumber\\
&\forall k\in\N_{N-1},\quad  \forall i\in\mc{I}_\mc{Y},\quad  c_i^Tz_k+d_i^Tv_k\le 1-f_i,\nonumber\\
&\forall k\in\N_{N_\mc{Z}-1},\quad  z_{N+k+1}=(A+BK_\mc{Z})z_{N+k},\text{ and}\nonumber\\
\label{eq:04.03}
&\forall k\in\N_{N_\mc{Z}},\quad  \forall i\in\mc{I}_\mc{Y},\quad  (c_i^T+d_i^TK_\mc{Z})z_{N+k}\le 1-f_i.
\end{align}
It is worth observing that the numbers of affine equalities and inequalities in the constraints of~\eqref{eq:04.03} is
\begin{align}
\label{eq:04.04}
n_{\mbf{eq}_T}&=n+Nn+N_\mc{Z}n,\text{ and}\nonumber\\
n_{\mbf{iq}_T}&=N_\mc{S}q+Np+(N_\mc{Z}+1)p.
\end{align}
In what follows, $V_T^0\bb$, $\mbf{d}_T^0\bb$ and $\mc{C}_T$  denote, respectively, the value function and the optimizer map of the optimization problem~\eqref{eq:04.03}, and their effective domain.  The set $\mc{C}_T$, also referred to as the $T$--controllable set, is the set of states $x$ for which there exists at least one decision variable $\mbf{d}_T$ satisfying constraints in~\eqref{eq:04.03}. In our setting, $\mc{C}_T$ is a closed polyhedral proper $D$--set in $\R^n$, and its interior contains the set $\mc{S}$. Likewise, in our setting, the cost function $\mbf{d}_T\mapsto V_N(\mbf{d}_T)$ is a convex and quadratic function, which is strictly convex with respect to the $(\mbf{z}_{N-1},\mbf{v}_{N-1},\mbf{z}_{[N:N+N_\mc{Z}]})$ component of the decision variable $\mbf{d}_T$. Furthermore, by the definitions of the stage and overall cost functions, there exists scalars $\beta_1\in (0,\infty)$ and $\beta_3\in (0,\infty)$ such that, for all $(z,v)\in\R^n\times \R^m$, and, for all $\mbf{d}_T\in\R^{n_{\mbf{d}_T}}$, with $z=z_0$,
\begin{align*}
&\beta_1\|z\|^2\le \ell(z,v),\\
&\beta_1\|z\|^2\le \beta_1\|(\mbf{z}_{N-1},\mbf{v}_{N-1},\mbf{z}_{[N:N+N_\mc{Z}]})\|^2,\\
&\beta_1\|(\mbf{z}_{N-1},\mbf{v}_{N-1},\mbf{z}_{[N:N+N_\mc{Z}]})\|^2\le V_T(\mbf{d}_T)\ \text{and}\\
&V_T(\mbf{d}_T)\le \beta_3\|(\mbf{z}_{N-1},\mbf{v}_{N-1},\mbf{z}_{[N:N+N_\mc{Z}]})\|^2.
\end{align*}
Thus, the implicit rigid tube optimal control problem $\mathfrak{P}_T(x)$ is a convex quadratic programming problem for any given $x\in\R^n$, which is feasible for any $x\in\mc{C}_T$. In light of the forms of the cost function~\eqref{eq:04.02} and the constraints in~\eqref{eq:04.03}, this convex quadratic programming problem is highly structured, i.e., the cost and equality and inequality constraints in its formulation can be expressed by using very sparse and almost block diagonal matrices. In fact, the optimization problem $\mathfrak{P}_T(x)$ is a parametric convex quadratic programming problem with respect to the state $x$, and the structure of its solution is well understood and documented~\cite{rockafellar:wets:2009}. 
\begin{theorem}
\label{thm:04.01}
Suppose Assumptions~\ref{ass:02.01}--\ref{ass:02.10} hold, and take any $N\ge 1$. 
\begin{itemize}
\item [$(i)$] The optimization problem $\mathfrak{P}_T(x),\ x\in\R^n$ is a convex quadratic programming problem, which is feasible for all $x\in\mc{C}_T$.
\item [$(ii)$] The value function $V_T^0\bb \ :\ \mc{C}_T\rightarrow \R$ is a convex, continuous and piecewise quadratic function such that 
\begin{align*}
\forall x\in\mc{S},&\  V_T^0(x)=0\ \text{and}\\
\forall x\in\mc{C}_T\setminus \mc{S},&\  0<V_T^0(x)<\infty.
\end{align*}
\item [$(iii)$] The components $\mbf{z}_{N-1}^0\bb\ :\ \mc{C}_T\rightarrow \R^{Nn}$, $\mbf{v}_{N-1}^0\bb\ :\ \mc{C}_T\rightarrow \R^{Nm}$ and  $\mbf{z}_{[N:N+N_\mc{Z}]}^0\bb\ :\ \mc{C}_T\rightarrow \R^{(N_\mc{Z}+1)n}$ of the optimizer map $\mbf{d}_T^0\bb\ :\  \mc{C}_T\rightarrow \R^{n_{\mbf{d}_T}}$ are continuous and piecewise affine functions.
\item [$(iv)$] There exists a continuous and piecewise affine function $\boldsymbol{\omega}_{N_\mc{S}-1}^*\bb\ :\ \mc{C}_T\rightarrow \R^{N_\mc{S}n}$ such that, for all $x\in\mc{C}_T$,
\begin{equation*}
(\mbf{z}_{N-1}^0,\mbf{v}_{N-1}^0,\mbf{z}_{[N:N+N_\mc{Z}]}^0,\boldsymbol{\omega}_{N_\mc{S}-1}^*)(x)\in \mbf{d}_T^0(x).
\end{equation*}
\end{itemize}
\end{theorem}

It is worth noting that the functions $\mbf{z}_{N-1}^0\bb$, $\mbf{v}_{N-1}^0\bb$ and  $\mbf{z}_{[N:N+N_\mc{Z}]}^0\bb$ are, in fact, Lipschitz continuous since they are continuous and piecewise affine functions defined over a nonempty closed polyhedral set $\mc{C}_T$.
\begin{corollary}
\label{cor:04.01}
Suppose Assumptions~\ref{ass:02.01}--\ref{ass:02.10} hold, and take any $N\ge 1$.   The components $\mbf{z}_{N-1}^0\bb\ :\ \mc{C}_T\rightarrow \R^{Nn}$, $\mbf{v}_{N-1}^0\bb\ :\ \mc{C}_T\rightarrow \R^{Nm}$ and  $\mbf{z}_{[N:N+N_\mc{Z}]}^0\bb\ :\ \mc{C}_T\rightarrow \R^{(N_\mc{Z}+1)n}$ of the optimizer map $\mbf{d}_T^0\bb\ :\  \mc{C}_T\rightarrow \R^{n_{\mbf{d}_T}}$ are, in fact, Lipschitz continuous  functions.
\end{corollary}

\subsection{Implicit Rigid Tube MPC}
\label{sec:04.03}
The implicit rigid tube MPC law $\kappa_T\bb\ :\ \mc{C}_T\rightarrow \R^m$ is specified, for all $x\in\mc{C}_T$, by
\begin{equation}
\label{eq:04.08}
\kappa_T(x)=v_0^0(x)+K_\mc{S}(x-z_0^0(x)),
\end{equation}
and the related implicit rigid tube model predictive controlled dynamics is given, for all $x\in\mc{C}_T$, by
\begin{equation}
\label{eq:04.09}
x^+\in Ax+B\kappa_T(x)\oplus \mc{W}.
\end{equation}
By Theorem~\ref{thm:04.01} and Corollary~\ref{cor:04.01}, the control function $x\mapsto \kappa_T(x)$ is piecewise affine and Lipschitz continuous relative to the $T$--controllable set $\mc{C}_T$. Likewise, the closed--loop dynamics map $x\mapsto Ax+B\kappa_T(x)\oplus \mc{W}$ is a polyhedral, polytopic--valued and Lipschitz continuous, with respect to the Hausdorff distance, set--valued map relative to the $T$--controllable set $\mc{C}_T$.  In the actual implementation, $v_0^0(x)$ and $z_0^0(x)$ are the first terms of the sequences $\mbf{z}_{N-1}^0(x)=\{z_k^0(x)\}_{k=0}^{N-1}$ and $\mbf{v}_{N-1}^0(x)=\{v_k^0(x)\}_{k=0}^{N-1}$, which are computed by solving online the implicit rigid tube optimal control problem~\eqref{eq:04.03}.

The robust positive invariance and exponential stability properties of the implicit rigid tube MPC are inherited from its predecessor~\cite{mayne:seron:rakovic:2005}, as briefly summarized next.  We recall that a (nonempty) set $\mc{C}$ is a robust positively invariant set for the uncertain dynamics $x^+\in f(x)\oplus \mc{W}$ and constraints $x\in\mc{X}$ if and only if for all $x\in\mc{C}$ it holds that $f(x)\oplus \mc{W}\subseteq \mc{C}$ and $x\in\mc{X}$. In our setting, the $T$--controllable set $\mc{C}_T$ is a robust positively invariant set. 
\begin{theorem}
\label{thm:04.02}
Suppose Assumptions~\ref{ass:02.01}--\ref{ass:02.10} hold, and take any $N\ge 1$. The $T$--controllable set $\mc{C}_T$ is a robust positively invariant set for the implicit rigid tube model predictive controlled uncertain dynamics  $x^+\in Ax+B\kappa_T(x)\oplus \mc{W}$ and constraints $(x,\kappa_T(x))\in\mc{Y}$, i.e.,
\begin{equation*}
\forall x\in\mc{C}_T,\  Ax+B\kappa_T(x)\oplus \mc{W}\subseteq \mc{C}_T\  \text{and}\  (x,\kappa_T(x))\in\mc{Y}.
\end{equation*}
\end{theorem}

We recall that the point--to--set distance function is defined, for any set $\mc{S}\subseteq \R^n$ and any point $x\in\R^n$, by 
\begin{equation*}
\operatorname{distance}(\mc{S},x):=\inf_{s}\{\|x-s\|\ :\ s\in\mc{S}\}.
\end{equation*}
 We also recall that when (nonempty) closed sets $\mc{S}$ and $\mc{C}$ in $\R^n$ such that $\mc{S}\subset \operatorname{interior}(\mc{C})$ are robust positively invariant  for the uncertain dynamics $x^+\in f(x)\oplus \mc{W}$ and constraints $x\in\mc{X}$ and when there exist a function $V\bb\ :\ \mc{C}\rightarrow \R$ and two scalars $\gamma_1\in(0,\infty)$ and $\gamma_2\in(0,\infty)$ such that
\begin{align*}
\forall x\in \mc{C},\  &\gamma_1\operatorname{distance}(\mc{S},x)\le V(x)\le \gamma_2\operatorname{distance}(\mc{S},x),\text{ and}\\
\forall x\in \mc{C},\  &\forall x^+\in f(x)\oplus \mc{W},\\
 &V(x^+)\le V(x) -\gamma_1\operatorname{distance}(\mc{S},x),
\end{align*}
then the set $\mc{S}$ is robustly exponentially stable  for the uncertain dynamics $x^+\in f(x)\oplus \mc{W}$  with the domain of attraction equal to the set $\mc{C}$.  In our setting, the previously outlined properties of the overall cost function $V_T\bb$, Theorem~\ref{thm:04.01} and Corollary~\ref{cor:04.01} yield the following.
\begin{proposition}
\label{prop:04.03}
Suppose Assumptions~\ref{ass:02.01}--\ref{ass:02.10} hold, and take any $N\ge 1$. There exists scalars $\beta_1\in (0,\infty)$ and $\beta_2\in(0,\infty)$ such that, for all $x\in\mc{C}_T$,
\begin{equation*}
\beta_1\|z_0^0(x)\|^2\le V_T^0(x)\le \beta_2\|z_0^0(x)\|^2,
\end{equation*}
and, for all $x\in\mc{C}_T$ and for all $x^+\in Ax+B\kappa_N(x)\oplus \mc{W}$,
\begin{equation*}
V_T^0(x^+)\le V_T^0(x) -\beta_1\|z_0^0(x)\|^2.
\end{equation*}
\end{proposition}
The preceding facts ensure that the implicit rigid tube MPC inherits the desirable robust exponential stability of the set $\mc{S}$ from the ordinary rigid tube MPC~\cite{mayne:seron:rakovic:2005}.
\begin{theorem}
\label{thm:04.03}
Suppose Assumptions~\ref{ass:02.01}--\ref{ass:02.10} hold. Take any $N\ge 1$. The set $\mc{S}$ is robustly exponentially stable for the implicit rigid tube model predictive controlled uncertain dynamics  $x^+\in Ax+B\kappa_T(x)\oplus \mc{W}$ with the domain of attraction equal to the $T$--controllable set $\mc{C}_T$.
\end{theorem}

\section{A Concise Step--By--Step Design Guidelines}
\label{sec:05}

\subsection{The Matrix $K_\mc{S}$}
\label{sec:05.01}

Since, by Assumption~\ref{ass:02.01}, the matrix pair $(A,B)$ is strictly stabilizable, it is always possible to select a control matrix $K_\mc{S}\in\Rnm{m}{n}$ such that the matrix $A+BK_\mc{S}$ is strictly stable, as postulated in Assumption~\ref{ass:02.05}. 

\subsection{The Set $\mc{S}$}
\label{sec:05.02}

The implicit representation~\eqref{eq:02.10} of  the set $\mc{S}$ postulated in Assumption~\ref{ass:02.06} requires only a scalar $\alpha\in [0,1)$ and a finite integer $N_\mc{S}$, for which the condition~\eqref{eq:02.11} holds true, to be determined. This can be done without the need to construct explicitly the sets $(A+BK_\mc{S})^{N_\mc{S}}\mc{W}$. Indeed, in light of Assumption~\ref{ass:02.03}, for any  $\alpha\in [0,1)$ and any $N_\mc{S}\in \N$, the set inclusion 
\begin{equation*}
(A+BK_\mc{S})^{N_\mc{S}}\mc{W}\subseteq \alpha \mc{W},
\end{equation*}
holds true if and only if  
\begin{align}
\label{eq:05.01}
\alpha_{N_\mc{S}}&:=\max_{i\in\mc{I}_\mc{W}} \alpha_{(N_\mc{S},i)} \le \alpha,\text{ with}\\
\label{eq:05.02}
\forall i\in\mc{I}_\mc{W},\  \alpha_{(N_\mc{S},i)}&:=\support{\mc{W}}{((A+BK_\mc{S})^{N_\mc{S}})^\top e_i}.
\end{align}
Since $0\in\operatorname{interior}(\mc{W})$, $\alpha_{(N_\mc{S},i)}\ge 0$ for all $i\in\mc{I}_\mc{W}$. Thus, the condition~\eqref{eq:02.11} holds true if and only if $N_\mc{S}\in\N$ and $\alpha_{N_\mc{S}}\le\alpha\in[0,1)$, which, for each fixed $N_\mc{S}\in\N$ and any given $\alpha \in [0,1)$,  can be checked by making use of:
\begin{itemize}
\item the solution to $q$ simple linear programming problems,  which are used to compute the corresponding values $\support{\mc{W}}{((A+BK_\mc{S})^{N_\mc{S}})^\top e_i},\ i\in\mc{I}_\mc{W}$ in order to evaluate the scalars $\alpha_{(N_\mc{S},i)},\ i\in\mc{I}_\mc{W}$ in~\eqref{eq:05.02},
\item  the computation of the maximum $\alpha_{N_\mc{S}}$ of these $q$ scalars $\alpha_{(N_\mc{S},i)},\ i\in\mc{I}_\mc{W}$, and
\item  the comparison of the maximum $\alpha_{N_\mc{S}}$ to $\alpha$.
\end{itemize}
As for the linear programming problems needed to evaluate $\support{\mc{W}}{((A+BK_\mc{S})^{N_\mc{S}})^\top e_i},\ i\in\mc{I}_\mc{W}$  in~\eqref{eq:05.02},
\begin{itemize}
\item the total number of decision variables is $n$,
\item the total number of affine equalities is $0$, and 
\item the total number of affine inequalities is $q$.
\end{itemize}
In practical terms, one usually fixes a value of $\alpha\in (0,1)$ and identifies the smallest integer $N_\mc{S}$ for which~\eqref{eq:05.01} holds true, and one utilizes the corresponding value of $\alpha_{N_\mc{S}}$ as an optimal value for $\alpha$, i.e., one sets $\alpha=\alpha_{N_\mc{S}}$.

The verification of Assumption~\ref{ass:02.07} can be also done without the need to perform the explicit computation of the set $\mc{S}$ specified in~\eqref{eq:02.10}. Namely, in light of relations in~\eqref{eq:03.05}, the scalars $f_{(N_\mc{S},i)}:=\support{\mc{S}}{c_i+K_\mc{S}^\top d_i},\ i\in\mc{I}_\mc{Y}$ can be efficiently computed by
\begin{equation}
\label{eq:05.03}
f_{(N_\mc{S},i)}=(1-\alpha)^{-1}\sum_{j=0}^{N_\mc{S}-1}\support{\mc{W}}{((A+BK_\mc{S})^j)^\top \eta_i},
\end{equation}
where, for all $i\in\mc{I}_\mc{Y}$, $\eta_i:=c_i+K_\mc{S}^\top d_i$. The condition~\eqref{eq:02.12} holds true if and only if
\begin{equation}
\label{eq:05.04}
f_{N_\mc{S}}:=\max_{i\in\mc{I}_\mc{Y}}f_{(N_\mc{S},i)}<1.
\end{equation}
Since $0\in\operatorname{interior}(\mc{S})$, $f_{(N_\mc{S},i)}\ge 0$ for all  $i\in\mc{I}_\mc{Y}$. Consequently, for each fixed $N_\mc{S}\in\N$ and fixed $\alpha \in [0,1)$,  the condition~\eqref{eq:02.12} can be checked by making use of: 
\begin{itemize}
\item the solution to $N_\mc{S}p$ simple linear programming problems, which are used to compute the corresponding values $\support{\mc{W}}{((A+BK_\mc{S})^j)^\top \eta_i},\ j\in N_{N_\mc{S}-1},\ i\in\mc{I}_\mc{Y}$ in order to evaluate the scalars $f_{(N_\mc{S},i)},\ i\in\mc{I}_\mc{Y}$ in~\eqref{eq:05.03},
\item  the computation of the maximum $f_{N_\mc{S}}$ of these $p$ scalars $f_{(N_\mc{S},i)},\ i\in\mc{I}_\mc{Y}$, and
\item  the comparison of the maximum $f_{N_\mc{S}}$ to $1$.
\end{itemize}

Clearly, the scalars $f_{(N_\mc{S},i)}=\support{\mc{S}}{c_i+K_\mc{S}^\top d_i},\ i\in\mc{I}_\mc{Y}$ can be also computed by a simple iteration given, for all $k\in\N$ and all $i\in\mc{I}_\mc{Y}$, with $f_{(0,i)}=0$, by
\begin{equation}
\label{eq:05.05}
f_{(k+1,i)}=f_{(k,i)}+(1-\alpha)^{-1}\support{\mc{W}}{((A+BK_\mc{S})^k)^\top \eta_i}.
\end{equation}
Thus, in practical terms, the conditions~\eqref{eq:05.01} and~\eqref{eq:05.04} can be verified either independently or jointly via an adequately defined iteration combining~\eqref{eq:05.01} and~\eqref{eq:05.05}. 

As for the linear programming problems needed to evaluate either $\support{\mc{W}}{((A+BK_\mc{S})^j)^\top \eta_i},\ i\in\mc{I}_\mc{Y}$ in~\eqref{eq:05.03} or $\support{\mc{W}}{((A+BK_\mc{S})^k)^\top \eta_i},\ i\in\mc{I}_\mc{Y}$ in~\eqref{eq:05.05},
\begin{itemize}
\item the total number of decision variables is $n$,
\item the total number of affine equalities is $0$, and 
\item the total number of affine inequalities is $q$.
\end{itemize}

Finally, if one wishes to deploy the set $\mc{S}$ that is a $\varepsilon$--robust positively invariant approximation~\cite{rakovic:kerrigan:kouramas:mayne:2004b} of the minimal robust positively invariant set, the condition~\eqref{eq:02.11} should be accompanied with the condition that $\alpha \mc{S}\subseteq \varepsilon \mc{B}$, where $\mc{B}$ is a closed unit norm ball. In the case of the polytopic norms, these two conditions can be also handled directly by  a suitably defined iteration combining~\eqref{eq:05.01} and a direct modification of~\eqref{eq:05.05}.

\subsection{The Set $\mc{Y}_\mc{S}$}
\label{sec:05.03}
The possibly redundant representation of the polyhedral proper $D$--set $\mc{Y}_\mc{S}$ is given directly by~\eqref{eq:03.04}. 

\subsection{The Matrices $K_\mc{Z}$ and $P$}
\label{sec:05.04}

Since, by Assumption~\ref{ass:02.01}, the matrix pair $(A,B)$ is strictly stabilizable, it is always possible to select a control matrix $K_\mc{Z}\in\Rnm{m}{n}$ such that the matrix $A+BK_\mc{Z}$ is strictly stable, as postulated in Assumption~\ref{ass:02.08}. An optimal design selection for the matrix $K_\mc{Z}\in\Rnm{m}{n}$ and the terminal cost weighting matrix $P$ is given from the solution of the unconstrained infinite horizon discrete time linear quadratic regulator for the disturbance free system $z^+=Az+Bv$ with the stage cost $\ell(z,v)=z^\top Q z+v^\top R v$. Alternatively, the control matrix $K_\mc{Z}\in\Rnm{m}{n}$ and the terminal cost weighting matrix $P$ can be chosen with ease to satisfy condition~\eqref{eq:02.18} by using numerous standard control techniques.

\subsection{The Set $\mc{Z}_\mc{S}$}
\label{sec:05.05}
For any control matrix $K_\mc{Z}\in\Rnm{m}{n}$ satisfying Assumptions~\ref{ass:02.08}, the possibly redundant representation of the polyhedral proper $D$--set $\mc{Z}_\mc{S}$ is given directly by~\eqref{eq:03.07}.

\subsection{The Set $\mc{Z}_f$}
\label{sec:05.06}

The implicit representation~\eqref{eq:02.16} of the terminal constraint set $\mc{Z}_f$ postulated in Assumption~\ref{ass:02.09} requires only a finite integer $N_\mc{Z}$ to be determined for which the condition~\eqref{eq:02.17} holds. This condition read as 
\begin{equation*}
\bigcap_{j=0}^{N_\mc{Z}}(A+BK_\mc{Z})^{-j}\mc{Z}_\mc{S}\subseteq (A+BK_\mc{Z})^{-(N_\mc{Z}+1)}\mc{Z}_\mc{S}
\end{equation*}
and, thus, for a fixed $N_\mc{Z}\in\N$, it holds true if and only if
\begin{equation}
\label{eq:05.06}
\forall i\in\mc{I}_\mc{Y},\  \support{\bigcap_{j=0}^{N_\mc{Z}}(A+BK_\mc{Z})^{-j}\mc{Z}_\mc{S}}{\psi_i}\le 1-f_i,
\end{equation}
where each $\psi_i:=((A+BK_\mc{Z})^{(N_\mc{Z}+1)})^\top (c_i + K_\mc{Z}^\top d_i)$, and which requires the evaluation of the support function of a set given as the set intersection of the finitely many preimage sets. Such an evaluation of the support function can be performed efficiently, and it does not require the considered set to be constructed explicitly. 
\begin{lemma}
\label{lemma:05.01}
Let $\mc{X}$ and $M$ be a proper $D$--set in $\R^n$ and a matrix in $\Rnm{n}{n}$, respectively. Let also $\mc{J}:=\{0,1,\ldots,j\}$ for a finite $j\in\N,\ j>0$.  For all $y\in\R^n$, the value $\support{\bigcap_{j\in\mc{J}} M^{-j}\mc{X}}{y}$ of the support function $\supportbb{\bigcap_{j\in\mc{J}} M^{-j}\mc{X}}$ evaluated at $y$ is the maximum value of the following convex optimization problem
\begin{align*}
\text{maximize}\ &y^\top x\nonumber\\
\text{with respect to}\ &(x,\{x_j\}_{j\in\mc{J}})\nonumber\\
\text{subject to}\ &\forall j\in\mc{J},\ x_j=M^jx,\text{ and}\nonumber\\
&\forall j\in\mc{J},\ x_j\in\mc{X}.
\end{align*}
\end{lemma}
In light of Lemma~\ref{lemma:05.01}, when $N_{\mc{Z}}\ge 1$, the value $\support{\bigcap_{j=0}^{N_\mc{Z}}(A+BK_\mc{Z})^{-j}\mc{Z}_\mc{S}}{\psi}$ of the support function $\supportbb{\bigcap_{j=0}^{N_\mc{Z}}(A+BK_\mc{Z})^{-j}\mc{Z}_\mc{S}}$ is given, for all $\psi\in\R^n$, by the solution to the linear programming problem
\begin{align}
&\text{maximize}\ \psi^\top z_0\nonumber\\
&\text{with respect to}\ \{z_k\in\R^n\}_{k=0}^{N_{\mc{Z}}}\nonumber\\
&\text{subject to}\nonumber\\&\forall k\in\N_{N_{\mc{Z}}-1},\  z_{k+1}=(A+BK_\mc{Z})z_k,\text{ and}\nonumber\\
\label{eq:05.07}
&\forall k\in\N_{N_{\mc{Z}}},\  \forall i\in\mc{I}_\mc{Y},\   (c_i + K_\mc{Z}^\top d_i) z_k\le 1-f_i.
\end{align}
When $N_\mc{Z}=0$, $\bigcap_{j=0}^{N_\mc{Z}}(A+BK_\mc{Z})^{-j}\mc{Z}_\mc{S}=\mc{Z}_S$ and, for all $\psi\in\R^n$, $\support{\bigcap_{j=0}^{N_\mc{Z}}(A+BK_\mc{Z})^{-j}\mc{Z}_\mc{S}}{\psi}=\support{\mc{Z}_\mc{S}}{\psi}$. With this in mind, let,  for any $N_\mc{Z}\in\N$ and all $i\in\mc{I}_\mc{Y}$,
\begin{equation}
\label{eq:05.08}
g_{(N_\mc{Z},i)}:=\support{\bigcap_{j=0}^{N_\mc{Z}}(A+BK_\mc{Z})^{-j}\mc{Z}_\mc{S}}{\psi_i},
\end{equation}
where each $\psi_i=((A+BK_\mc{Z})^{(N_\mc{Z}+1)})^\top (c_i + K_\mc{Z}^\top d_i)$.  It follows that, for any given $N_\mc{Z}\in\N$, the condition~\eqref{eq:02.17} holds true if and only if
\begin{equation}
\label{eq:05.09}
\forall i\in\mc{I}_\mc{Y},\  g_{(N_\mc{Z},i)}+f_i\le 1,
\end{equation}
and that it can be checked efficiently via standard linear programming.
Summarizing, for each fixed $N_\mc{Z}\in\N$, the condition~\eqref{eq:02.17} can be checked by making use of:
\begin{itemize}
\item the solution to $p$ linear programming problems,  which are used to compute the corresponding values $\support{\bigcap_{j=0}^{N_\mc{Z}}(A+BK_\mc{Z})^{-j}\mc{Z}_\mc{S}}{\psi_i},\ i\in\mc{I}_\mc{Y}$ when $N_\mc{Z}\ge 1$ and the values $\support{\mc{Z}_\mc{S}}{\psi_i},\ i\in\mc{I}_\mc{Y}$ when $N_\mc{Z}=0$ in order to evaluate the scalars $g_{(N_\mc{Z},i)},\ i\in\mc{I}_\mc{Y}$ in~\eqref{eq:05.08},
\item  the computation of the sums $g_{(N_\mc{Z},i)}+f_i,\ i\in\mc{I}_\mc{Y}$ of these $p$ scalars $g_{(N_\mc{Z},i)},\ i\in\mc{I}_\mc{Y}$ and previously computed scalars $f_i,\ i\in\mc{I}_\mc{Y}$, and
\item  the comparisons of the  sums $g_{(N_\mc{Z},i)}+f_i,\ i\in\mc{I}_\mc{Y}$ to $1$.
\end{itemize}
As for the linear programming problems needed to evaluate $\support{\bigcap_{j=0}^{N_\mc{Z}}(A+BK_\mc{Z})^{-j}\mc{Z}_\mc{S}}{\psi_i},\ i\in\mc{I}_\mc{Y}$ in~\eqref{eq:05.08},
\begin{itemize}
\item the total number of decision variables is $(N_{\mc{Z}}+1)n$,
\item the total number of affine equalities is $N_{\mc{Z}}n$, and 
\item the total number of affine inequalities is $(N_{\mc{Z}}+1)p$.
\end{itemize}

A  numerically more convenient, sufficient condition~\cite{rakovic:zhang:2023:a} ensuring that~\eqref{eq:02.17} holds is given by
\begin{equation}
\label{eq:05.10}
\mc{Z}_\mc{S}\subseteq (A+BK_\mc{Z})^{-(N_\mc{Z}+1)}\mc{Z}_\mc{S}.
\end{equation}
This sufficient condition~\eqref{eq:05.10} might increase the required value for $N_\mc{Z}$ relative to the value required by the necessary and sufficient condition~\eqref{eq:02.17}, but its use can be  justified by its simplicity. Indeed, its is equivalent to
\begin{equation}
\label{eq:05.11}
\forall i\in\mc{I}_\mc{Y},\  \support{\mc{Z}_\mc{S}}{\psi_i}+f_i\le 1,
\end{equation}
where each $\psi_i:=((A+BK_\mc{Z})^{(N_\mc{Z}+1)})^\top (c_i + K_\mc{Z}^\top d_i)$, and, for each fixed $N_\mc{Z}\in\N$,  the condition~\eqref{eq:05.10} can be checked by making use of 
\begin{itemize}
\item the solution to $p$ linear programming problems,  which are used to compute the corresponding values $\support{\mc{Z}_\mc{S}}{\psi_i},\ i\in\mc{I}_\mc{Y}$,
\item  the computation of the sums $\support{\mc{Z}_\mc{S}}{\psi_i}+f_i,\ i\in\mc{I}_\mc{Y}$, and
\item  the comparisons of the computed sums $\support{\mc{Z}_\mc{S}}{\psi_i}+f_i,\ i\in\mc{I}_\mc{Y}$ to $1$.
\end{itemize}
As for the linear programming problems needed to evaluate $\support{\mc{Z}_\mc{S}}{\psi_i},\ i\in\mc{I}_\mc{Y}$ in~\eqref{eq:05.11},
\begin{itemize}
\item the total number of decision variables is $n$,
\item the total number of affine equalities is $0$, and 
\item the total number of affine inequalities is $p$.
\end{itemize}

\section{Discussion}
\label{sec:06}

\subsection{Numerical Experience}
\label{sec:06.01}
We highlight the efficiency of the implicit tube MPC by an overview of the offline computational effort needed for a sample of random control problems. For each pair in a collection of state and control  dimensions, a sample of $1000$ controllable control systems was generated using DRSS function in MATLAB. For convenience, for each example, a relatively simple stage constraint set $\mc{Y}=100\mc{B}_\infty^n\times 50 \mc{B}_\infty^m$ and a relatively simple disturbance bounding set $\mc{W}=\mc{B}_\infty^n$ were considered. In each example, the control matrix $K_\mc{S}$ was selected by placing the eigenvalues of the matrix $A+BK_\mc{S}$ to the locations specified by the following set $\{-\frac{1}{4}+\frac{i}{2(n-1)}\ :\ i\in\N_{n-1}\}$ so that  $\rho(A+BK_\mc{S})\le 0.25$. The maximal number of allowed iterations  used for the computations of $N_\mc{S}$, $\alpha$ and $\{f_i\}_{i=1}^p$ was set to $10000$ per example. Also, for simplicity, in each example $Q=I$ and $R=I$ were used, while $K_\mc{Z}$ and $P$ were determined as the solution of the unconstrained infinite horizon discrete time linear quadratic regulator for the disturbance free system $z^+=Az+Bv$ with the stage cost $\ell(z,v)=z^\top Q z+v^\top R v$.  Finally, in each example, the computationally simpler, but sufficient, condition~\eqref{eq:05.10} was utilized in order to identify a guaranteed value for $N_\mc{Z}$. 
\begin{table}[!hbt]    
\begin{center}
\resizebox{.4875\textwidth}{!}{%
  \begin{tabular}{c|ccccccccc}
  \hline
     $n$ &  $2$ & $3$& $5$ & $8$ &  $13$ & $21$ & $34$ & $55$ & $89$ \\ 
     $m$ &  $1$ & $1$& $1$ & $2$ &  $3$ & $5$ & $7$ &  $11$ & $18$ \\ 
    \hline
     $\widehat{N}_\mc{S}$ &  $2$ & $3$ & $4$ & $4$ & $5$ & $5$ & $6$ & $6$ & $6$\\  
     $\widehat{\alpha}$   &  $0.30$ & $0.42$ & $0.42$ & $0.38$ & $0.37$ & $0.36$ & $0.33$ & $0.30$ & $0.28$\\ 
     $\widehat{t}_\mc{S}$ &  $0.09$ & $0.08$ & $0.13$ & $0.17$ & $0.24$ & $0.37$ & $0.71$ & $1.3$ & $2.8$\\ 
     $\widehat{N}_\mc{Z}$ &  $2$ & $3$ & $5$ & $4$ & $5$ & $6$ & $8$ & $9$ & $10$\\
     $\widehat{t}_\mc{Z}$ &  $0.07$ & $0.14$ & $0.43$ & $0.678$ & $1.40$ & $2.44$ & $5.15$ & $10.05$ & $19.83$\\ 
     $sR$ &  $92\%$ & $81\%$ & $66\%$ & $82\%$ & $86\%$ & $92\%$ & $89\%$ & $91\%$ & $95\%$\\
      \hline    
  \end{tabular}
}
Table 1. Offline Computations for Sets of Random Systems.
\end{center}
\end{table}
Table 1 provides a summarized overview of the performed computations. In this table, $\widehat{N}_\mc{S}$ and $\widehat{\alpha}$ are (successful sample based) averages for $N_\mc{S}$ and $\alpha$ for each consider state and control dimension, while $\widehat{t}_\mc{S}$ is the average time in \emph{miliseconds} for the computation of the successful samples. Likewise, in this table, $\widehat{N}_\mc{Z}$ is (a successful sample based) average for $N_\mc{Z}$, and $\widehat{t}_\mc{Z}$ is average time in \emph{seconds} for the computation of the successful samples. Finally, $sR$ represents the success rate, which is a number of examples with successful completion of the offline design divided by total number ($1000$) of the considered random examples. The offline design stage fails either because the maximum number of allowed iterations is reached without completion of the offline design or  because the underlying Assumptions can not be satisfied for a specific control problem sample. Evidently, the data summarized in Table 1 provides very strong evidence of the computational efficiency of the developed implicit rigid tube MPC.  
\begin{table}[!hbt]    
\begin{center}
\resizebox{.4875\textwidth}{!}{%
  \begin{tabular}{c|ccccccccc}
      \hline    
     $n$ &   $144$ & & $233$ & & $377$ & & $610$ & & $987$\\ 
     $m$ &   $29$ & & $47$ & &  $76$ & & $122$ & & $198$ \\ 
    \hline    
     $N_\mc{S}$  & $7$ & & $5$ & & $7$ & & $7$ & & $7$\\ 
     $\alpha$    & $0.0393$ & & $0.4104$ & & $0.0485$ & & $0.0409$ & & $0.0372$\\ 
     $t_\mc{S}\ (ms)$  & $21.8$ & & $21.9$ & & $66.8$ & & $403.6$ & & $767.2$\\ 
     $N_\mc{Z}$  & $12$ & & $14$ & & $14$ & & $17$ & & $11$\\
     $t_\mc{Z}\ (s)$  & $48.5$ & & $148.6$ & & $529.0$ & & $3369.0$ & & $13200.0$\\ 
     \hline
  \end{tabular}
}
Table 2. Offline Computations for Single  Random Systems.
\end{center}
\end{table}
As documented in Table 2, the computational efficiency of our proposal is further supported by the data underpinning the offline design for single random control problems in very high dimensions. As for the data summarized in Table 2, each random control problem is considered within the setting of the above outlined numerical experiments. 

\subsection{Illustration}
\label{sec:06.02}
Our illustration is a sampled variation of a model for a modern transport airplane representing a benchmark example AC9 in an excellent collection of examples~\cite{leibfritz:2004}. In this  illustrative example, 
\begin{equation*}
(n,m)=(10,4)
\end{equation*}
and the matrices $A\in\Rnm{10}{10}$ and $B\in\Rnm{10}{4}$ for the discrete time variation of the system are obtained from the matrices $A$ and $B$ for the original continuous time system specified in~\cite[Example AC9, Page 6]{leibfritz:2004}. The discrete time version of the system is obtained by using the standard Euler discretization with the sampling period $h=0.5$ \emph{seconds}. The stage state constraint set $\mc{Y}$ and disturbance bounding set $\mc{W}$ are given by
\begin{equation*}
\mc{Y}=500\mc{B}_\infty^{10}\times50\mc{B}_\infty^4 \text{ and }\mc{W}=\mc{B}_\infty^{10}.
\end{equation*}
The stage cost function $\ell\bbb$ is specified via matrices
\begin{equation*}
Q=100I\text{ and }R=I.
\end{equation*}

The offline stage of the design process is again undertaken by following steps specified in Section~\ref{sec:05}. In this example, the linear state feedback $x\mapsto K_\mc{S}x$ is, for diversity of the design illustration, obtained from the solution of the unconstrained infinite horizon discrete time linear quadratic regulator for the disturbance free system $s^+=As+B\upsilon$ with the stage cost $s^\top I s+\upsilon^\top I \upsilon$.  The first step of the design yields the implicit form of the rigid state tube cross--section shape set $\mc{S}$. In this example setting, the values of $\alpha$ and $N_\mc{S}$ are given by 
\begin{equation*}
\alpha=0.0408\text{ and }N_\mc{S}=24,
\end{equation*}
and they have been identified, together with the scalars $f_i=\support{\mc{S}}{c_i+K_\mc{S}^\top d_i},\ i\in\mc{I}_\mc{Y}$, in $t_\mc{S}=2.2$ \emph{milliseconds}. 

The terminal linear state feedback $z\mapsto K_\mc{Z}z$ and the terminal cost function $z\mapsto z^\top Pz$ are selected as the solution of the unconstrained infinite horizon discrete time linear quadratic regulator for the disturbance free system $z^+=Az+Bv$ with the stage cost $\ell(z,v)=z^\top Q z+v^\top R v$. In this illustrative example, we have intentionally selected $K_\mc{S}\neq K_\mc{Z}$ to highlight a design difference relative to the ordinary rigid tube MPC~\cite{mayne:seron:rakovic:2005} in which $K_\mc{S}=K_\mc{Z}$. The next step of the design is to determine the integer $N_\mc{Z}$ that yields the implicit form of the terminal set. In this example setting, $N_\mc{Z}=15$ is determined in $t_\mc{Z}=3.3$ \emph{seconds} by employing the condition~\eqref{eq:05.10}. Thus, the offline design is performed successfully in less than $4$ \emph{seconds}.

\begin{figure}[!h]
\includegraphics[width=0.4875\textwidth]{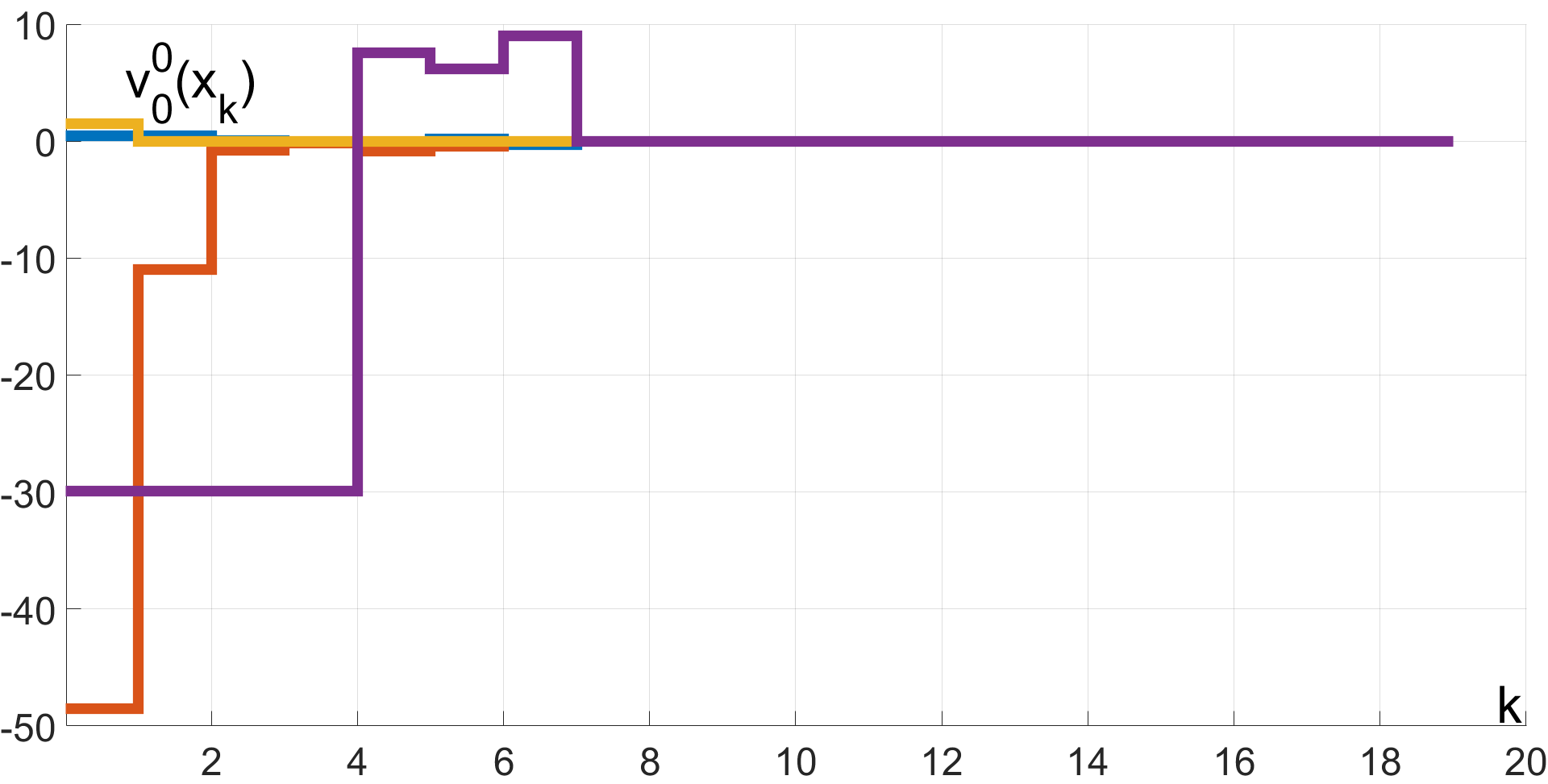}  
\caption{The $\{v_0^0(x_k)\}_{k=0}^{19}$ sequence.}
\label{fig:06.01}
\end{figure}

\begin{figure}[!h]
\includegraphics[width=0.4875\textwidth]{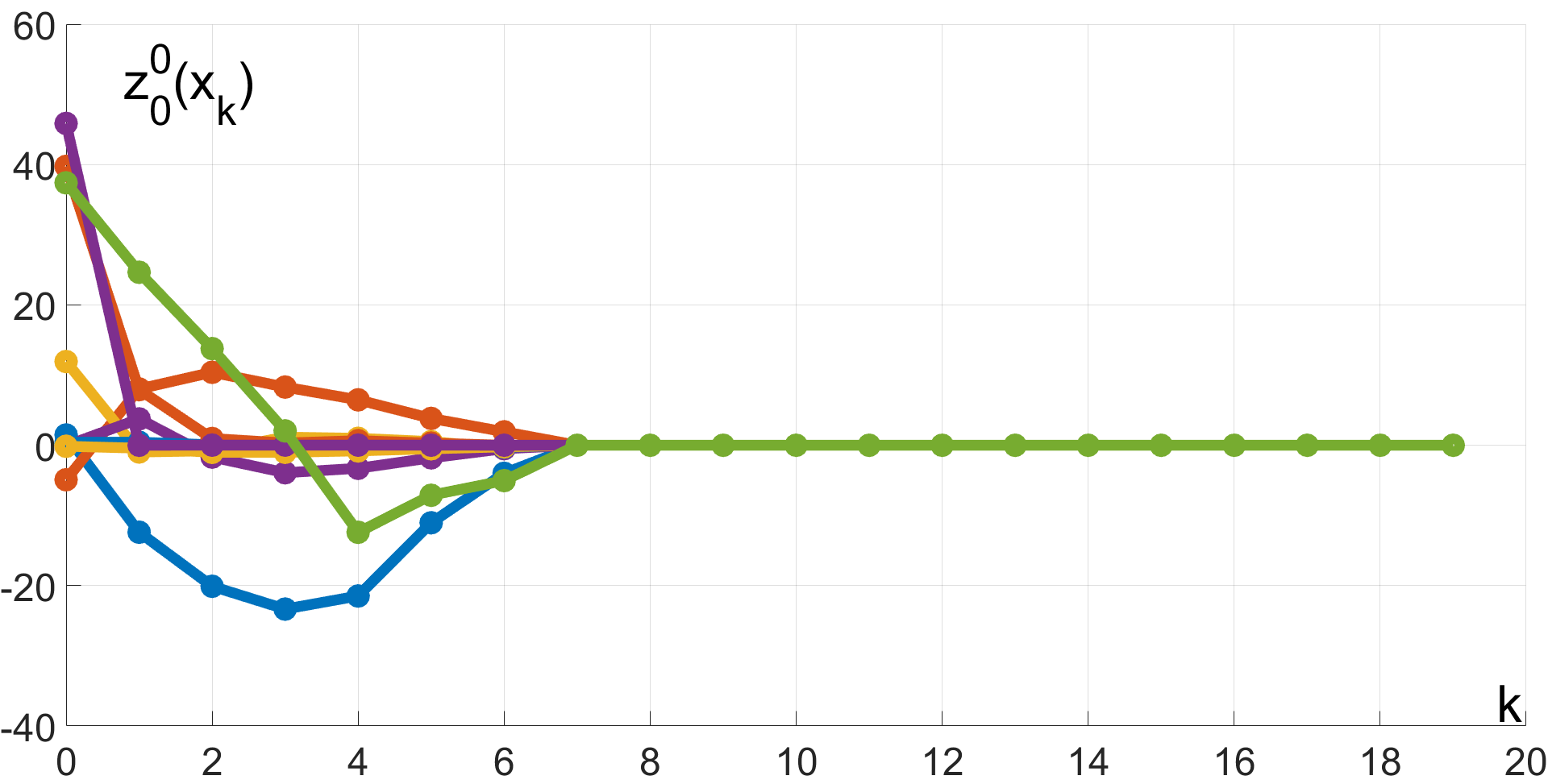}  
\caption{The $\{z_0^0(x_k)\}_{k=0}^{19}$ sequence.}
\label{fig:06.02}
\end{figure}

For the online implementation of the implicit rigid tube MPC, the prediction horizon length $N$ and the simulation time are both equal to $20$ time steps. The online optimization enabling the use of the implicit rigid tube model predictive takes the form of a convex quadratic programming problem, in which the total numbers of the decision variables, and equality and inequality constraints, are given, respectively, by
\begin{equation*}
n_{\mbf{d}_T}=680,\ n_{\mbf{eq}_T}=360\text{ and }n_{\mbf{iq}_T}=1488.
\end{equation*} 
This again illustrates the asserted computational efficiency of the implicit rigid tube MPC. 
The simulation results of the deployment of the implicit rigid tube MPC in a specific control process starting from the initial state 
\begin{align*}
x=(&7.3492, 24.3682, 20.3647, 40.7462, 5.0996,\\
 &10.5688, 42.1997, 0.1283, 46.8981, 29.7599)
\end{align*}
and affected by a sequence of randomly selected extreme disturbances are shown in Figures~\ref{fig:06.01}--\ref{fig:06.02}. 
 Figure~\ref{fig:06.01} depicts the sequence $\{v_0^0(x_k)\}_{k=0}^{19}$ and illustrates that the sequence $\{\|v_0^0(x_k)\|\}_{k=0}^{19}$ converges exponentially fast to $0$. This Figure indirectly illustrates that the values of the implicit rigid tube MPC $\kappa_T(x_k),\ k\in\N_{19}$ satisfy control constraints and converge exponentially fast to the set $K_\mc{S}\mc{S}$. 
Likewise,  Figure~\ref{fig:06.02} depicts the sequence $\{z_0^0(x_k)\}_{k=0}^{19}$ and illustrates that the sequence $\{\|z_0^0(x_k)\|\}_{k=0}^{19}$ converges exponentially fast to $0$. This Figure indirectly illustrates that the implicit rigid tube model predictive controlled states $x_k,\ k\in\N_{19}$ satisfy state constraints and converge exponentially fast to the set $\mc{S}$. We remark that in Figures~\ref{fig:06.01} and~\ref{fig:06.02}, only the sequences $\{v_0^0(x_k)\}_{k=0}^{19}$ and $\{z_0^0(x_k)\}_{k=0}^{19}$ are depicted in order to preserve clarity of these Figures. Once again, the observed behaviour of the implicit rigid tube MPC is, as expected, identical to that of the ordinary   rigid tube MPC~\cite{mayne:seron:rakovic:2005}. We close this illustration by pointing out that  the ordinary   rigid tube MPC~\cite{mayne:seron:rakovic:2005} fails to be applicable  on the account of the need to compute explicitly the sets $\mc{S}$ and $\mc{Z}_f$ and to use their explicit representations.

\subsection{Closing Remarks}
\label{sec:06.03}

This article has combined the classical works~\cite{gilbert:tan:1991,rakovic:kerrigan:kouramas:mayne:2004b,mayne:seron:rakovic:2005} on the maximal positively invariant set, the robust positively invariant approximations of the minimal robust positively invariant set and the rigid tube model predictive control with the more recent contributions~\cite{rakovic:2022,rakovic:zhang:2023:a,rakovic:zhang:2023:b} that have facilitated the use of the implicit set representations within the context of stability analysis, set invariance and MPC.  \emph{With regard to these and related works within the scope of the robust and tube MPC syntheses, the implicit set representations have not been previously utilized, and this article is, \emph{de facto}, the first work to pave the way for the use of these computationally efficient notions as well as document it properly}.

\emph{Summa summarum}, the computationally efficient implicit rigid tube MPC has been derived. The implicit set representations have been  deployed as a unique and novel feature of the implicit rigid tube MPC. This has been done in order to avoid the explicit computations of the rigid tubes and terminal constraint set deployed in the ordinary rigid tube MPC~\cite{mayne:seron:rakovic:2005}. As evident from a glance at the macroscopic computational summaries in terms of the structure and numbers of decision variables and constraints in Sections~\ref{sec:04} and~\ref{sec:05}, both the offline and online design stages of the implicit rigid tube MPC have been reduced to the standard, relatively simple, and reasonably structured linear and quadratic programming problems. In other words, relative to the ordinary rigid tube MPC~\cite{mayne:seron:rakovic:2005}, the implicit rigid tube MPC has been crafted in a manner that significantly improves the computational efficiency and retains scalability to high dimensions, effectively enables  a considerably broader practical utility, and preserves \emph{a priori} guarantees for all desirable topological and system theoretic properties.

An important line for future research is a generalization of the utilized implicit representations of the underlying tubes and terminal constraint sets within more general settings in terms of the considered control systems, constraints and costs.  Likewise, the development of the dedicated numerical tools for the underlying algebraic operations, and linear and quadratic programming problems, would further enhance the use of the implicit representations of the rigid tubes and terminal constraint sets. 

\textbf{Acknowledgement.} This is an extended version of the accepted Automatica brief paper $22$--$0760$, which was initially submitted on July $01$, $2022$ and was accepted for publication on June $22$, $2023$, and which is to be published by Automatica. The author is grateful to the Editor, Associate Editor and Referees for timely reviews and very constructive comments.

\section*{APPENDIX: Proofs and Evaluations of $\support{\mc{W}}{\eta}$ and $\support{\mc{Z}_\mc{S}}{\psi}$.}

\textbf{Proof of Lemma~\ref{lemma:03.01}.} Follows directly from the definitions of the Minkowski sum and set images.

\textbf{Proof of Lemma~\ref{lemma:03.02}.} The algebra of the support functions~\cite{rockafellar:1970,schneider:1993} yields, for all $y\in\R^n$,
\begin{equation*}
\support{\bigoplus_{j\in\mc{J}} M^j\mc{X}}{y}=\sum_{j\in\mc{J}}\support{M^j\mc{X}}{y}=\sum_{j\in\mc{J}}\support{\mc{X}}{(M^j)^\top y}.
\end{equation*}
\textbf{Proof of Lemma~\ref{lemma:03.03}.} For all $y\in\R^n$,  
$y\in \bigcap_{j\in\mc{J}} M^{-j}\mc{X}$ if and only if for all $j\in\mc{J},\  y\in M^{-j}\mc{X}$, i.e., if and only if for all $j\in\mc{J},\  M^jy\in \mc{X}$, which equivalently reads as  for all $j\in\mc{J},\  y_j\in \mc{X}\text{ with }y_j=M^jy$ and verifies the claim.

\textbf{Proof of Proposition~\ref{prop:03.01}}
In light of the definition of the sets $\mc{Y}_\mc{S}$ and $\mc{Y}$, $(z,v)\in\mc{Y}_\mc{S}$ if and only if for all $s\in \mc{S},\ (z+s,v+K_\mc{S}s)\in \mc{Y}$, which in terms of the support function $\supportbb{\mc{S}}$ reads equivalently as for all $i\in\mc{I}_\mc{Y},\ c_i^\top z+ d_i^\top v +\support{\mc{S}}{\eta_i}\le 1$ with $\eta_i:=c_i+K_\mc{S}^\top d_i$. Consequently, since, by Lemma~\ref{lemma:03.02}, for all $\eta \in\R^n$, 
\begin{equation*}
\support{\mc{S}}{\eta}=(1-\alpha)^{-1}\sum_{j=0}^{N_\mc{S}-1}\support{\mc{W}}{((A+BK_\mc{S})^j)^\top \eta},
\end{equation*}
it follows that
\begin{equation*}
\mc{Y}_\mc{S}=\{(z,v)\in \R^n\times \R^m\ :\ \forall i\in\mc{I}_\mc{Y},\  c_i^\top z + d_i^\top v \le 1-f_i\},
\end{equation*}
where, due to Assumption~\ref{ass:02.07}, for all $i\in\mc{I}_\mc{Y}$, the scalars
\begin{equation*}
f_i:=\support{\mc{S}}{\eta_i}=(1-\alpha)^{-1}\sum_{j=0}^{N_\mc{S}-1}\support{\mc{W}}{((A+BK_\mc{S})^j)^\top \eta_i}
\end{equation*}
are such that $f_i\in [0,1)$, so that $(1-f_i)\in (0,1]$ and, in turn, $\mc{Y}_\mc{S}$ is a polyhedral proper $D$--set in $\R^{n+m}$. 

\textbf{Validity of Theorem~\ref{thm:04.01}}
The claimed facts follow from the standard properties of the solution to the considered parametric convex quadratic programming problem and its construction. The interested reader is referred to~\cite{rockafellar:wets:2009} and~\cite{mayne:seron:rakovic:2005,rawlings:mayne:2009} for more details.

\textbf{Validity of Theorems~\ref{thm:04.02} and~\ref{thm:04.03} (and Proposition~\ref{prop:04.03})}
The claimed robust positive invariance and exponential stability properties of the implicit rigid tube MPC are both inherited from the ordinary rigid tube MPC~\cite{mayne:seron:rakovic:2005}.

\textbf{Proof of Lemma~\ref{lemma:05.01}.}  Since, for all $y\in\R^n$,
\begin{align*}
\support{\bigcap_{j\in\mc{J}} M^{-j}\mc{X}}{y}&=\sup_{x}\{y^\top x\ :\  x\in  \bigcap_{j\in\mc{J}} M^{-j}\mc{X}\},
\end{align*}
the claim follows from Lemma~\ref{lemma:03.03}.

\textbf{Evaluation of $\support{\mc{W}}{\eta}$.} In our setting, the evaluation of  $\support{\mc{W}}{\eta}$ can be performed by solving the following linear programming problem.
\begin{align*}
\text{maximize}\ &\eta^\top w\\
\text{with respect to}\ &w\\
\text{subject to}\ & \forall i\in\mc{I}_{\mc{W}},\ \e_i^\top w\le 1.
\end{align*}

\textbf{Evaluation of $\support{\mc{Z}_\mc{S}}{\psi}$.} In our setting, the evaluation of  $\support{\mc{Z}_\mc{S}}{\psi}$ can be performed by solving the following linear programming problem.
\begin{align*}
\text{maximize}\ &\psi^\top z\\
\text{with respect to}\ &z\\
\text{subject to}\ &  \forall i\in\mc{I}_\mc{Y},\  (c_i^\top + d_i^\top K_\mc{Z})z \le 1-f_i.
\end{align*}

\bibliographystyle{unsrt}
\bibliography{IRTMPC}
\end{document}